\documentclass[12pt]{amsart}
\usepackage{amscd}
\usepackage{verbatim}
\usepackage{amssymb, amsmath, amsthm, amscd,ifthen}
\usepackage[dvips]{graphics}
\usepackage[cp866]{inputenc}
\usepackage{graphicx}
\usepackage{epsfig}

\usepackage{amsmath,amssymb,amscd,}
\usepackage[all]{xy}

\usepackage{amssymb}

\newtheorem{thm}{Theorem}[section]

\newtheorem{prop}[thm]{Proposition}

\newtheorem{Ex}[thm]{Example}

\newtheorem{lemma}[thm]{Lemma}
\newtheorem{cor}[thm]{Corollary}

\theoremstyle{definition}

\newtheorem{dfn}[thm]{Definition}

\title
{Extremal  polygons in $\mathbb{R}^3$}

\author{ Gaiane Panina}

\address{
  St. Petersburg Institute for Informatics and Automation RAS;
St. Petersburg State University, e-mail:gaiane-panina@rambler.ru }


\begin{document}
\begin{abstract}
The oriented area function $A$ is (generically) a Morse function on
the
 space of  planar configurations of a polygonal linkage. We are lucky to
have an easy description of its critical points  as cyclic polygons
and a simple formula for the Morse index of a critical point.
However, for planar polygons, the function  $A$ in many cases is not a perfect
Morse function. In particular, for an equilateral pentagonal linkage
it has one extra local maximum (except for the global maximum) and
one extra local minimum.

In the present paper we consider the space of  3D configurations of
a polygonal linkage. For an appropriate generalization $S$ of the
area function $A$ the situation becomes
 nicer: we again have an easy description of  critical
points  and a simple formula for the Morse index.  In particular,
unlike the planar case, for an equilateral linkage with odd number
of edges the function $S$ is always a perfect Morse function and fits the
lacunary principle.  Therefore  cyclic equilateral polygons can be interpreted as
independent generators of the homology groups of the (decorated)
configuration space.
\end{abstract}

\maketitle \setcounter{section}{0}

\section{Introduction}
The oriented area function $A$ is (generically) a Morse function on
the
 space of  planar configurations of a polygonal linkage. We are lucky to
have an easy description of its critical points  as cyclic polygons (Theorem \ref{Thm_crirical_are_cyclic}), and a simple formula for
the Morse index of a critical point (Theorem
\ref{Thm_Morse_closed_plane}). However, for planar polygons, in many
cases $A$ is not a perfect Morse function. In particular, for an
equilateral pentagonal linkage it has one extra local maximum
(except for the global maximum) and one extra local minimum, see
Example \ref{Expentagon}.
  For an equilateral heptagonal linkage the number
of Morse points greatly exceeds the sum of Betti numbers of the
configuration space, and it is unclear how the boundary
homomorphisms of the Morse chain complex look like.

Surprisingly, if we pass to $\mathbb{R}^3$, for an appropriate
generalization $S$ of the area function $A$ the situation becomes
 nicer. We again have an easy description of  critical
points (as SW-invariant configurations, see Theorem
\ref{Thm_crirical_3D_SW}), and a simple formula for the Morse index
(Theorems \ref{Thm_Morse_planar}).  In particular, unlike the
planar case, for an equilateral linkage with odd number of edges
$S$, all critical points have even Morse indices. By the lacunary principle, $S$ is
 a perfect Morse function (see Theorem \ref{LemmaEquilateralTrue}) and
  the Morse chain complex  has zero boundary homomorphisms. As a
direct corollary we interpret cyclic equilateral polygons as
independent generators of  the homology groups of the configuration
space.

\textbf{Acknowledgements.} I'm indebted to George Khimshiashvili,
Dirk Siersma,  Alena Zhukova, and Mikhail Khristoforov for
inspiring conversations.

\section{Preliminaries and notation}\label{section_preliminaries}

A \textit{polygonal  $n$-linkage} is a sequence of positive numbers
$L=(l_1,\dots ,l_n)$. It should be interpreted as a collection of rigid
bars of lengths $l_i$ joined consecutively by revolving joints in a
chain.
\begin{dfn}
 \textit{A configuration} of $L$ in the
Euclidean space $ \mathbb{R}^d$, $d=2, \ 3$ is a sequence of points
$R=(p_1,\dots,p_{n+1}), \ p_i \in \mathbb{R}^d$ with
$l_i=|p_i,p_{i+1}|$ and $l_n=|p_n,p_{1}|$ modulo the action of
orientation preserving isometries of the space  $\mathbb{R}^d$.
We also call $P$ \textit{a closed chain} or a \textit{polygon}.

The set $M_d(L)$ of all  configurations  is \textit{the moduli
space, or the configuration space of the polygonal linkage }$L$.

\end{dfn}

 A configuration carries
a natural orientation which we indicate in figures by an arrow.

 We explain below in this paragraph what is known about planar
configurations and the signed area function as the Morse function on
the configuration space.

\begin{dfn} \label{Dfn_area} The \textit{signed area} of a polygon $P$ with the vertices \newline $p_i = (x_i,
y_i)$  is defined by
$$2A(P) = (x_1y_2 - x_2y_1) + \ldots + (x_ny_1 - x_1y_n).$$
\end{dfn}

\begin{dfn}
    A
polygon  $P$  is called \textit{cyclic} if all its vertices $p_i$
lie on a circle.

\end{dfn}

 Cyclic polygons  arise here as critical points of the signed area:

\begin{thm}\label{Thm_crirical_are_cyclic}\cite{khipan}
    Generically, a polygon $P$ is a critical point of the
signed area function $A$  iff $P$ is a cyclic configuration.
         \qed
\end{thm}

\begin{thm}\label{Thm_Morse_closed_plane}\cite{oberwolfach},\cite{zh}
For a generic cyclic configuration $P$ of a linkage $L$,
$$\mu(P)=\left\{
       \begin{array}{ll}
         e(P)-1-2\omega_P &\hbox{if }\   \delta(P)>0; \\
         e(P)-2-2\omega_P & \hbox{otherwise}.
       \end{array}
     \right.$$
     Here we used the below notation:
\end{thm}

\begin{figure}[h]
\centering
\includegraphics[width=6 cm]{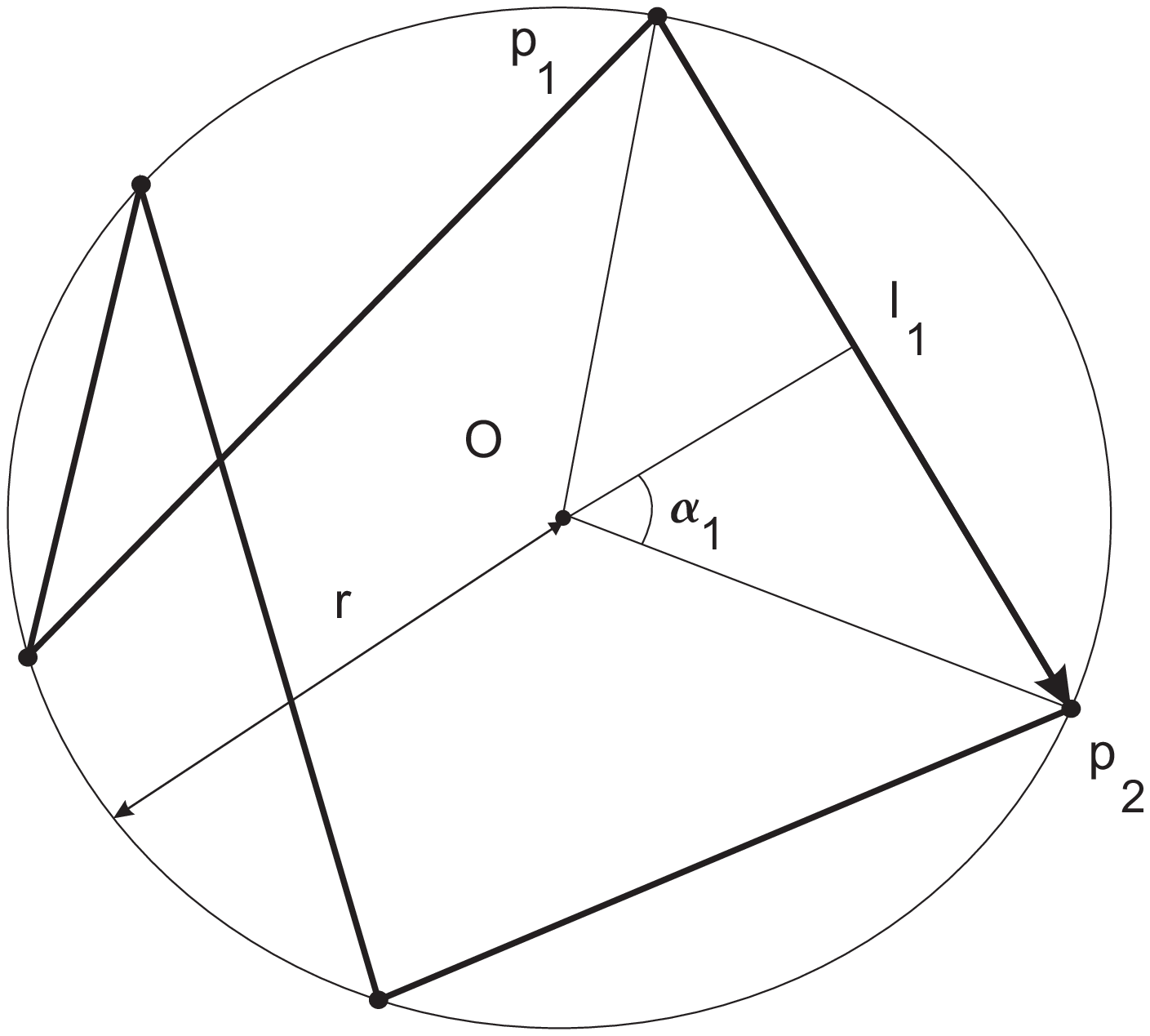}
\caption{Basic notation for a pentagonal cyclic configuration with
$E=(-1,-1,-1,1,-1)$}\label{Figure_notation}
\end{figure}

\textbf{Notation for cyclic
configurations, see Fig. \ref{Figure_notation}.}

$r$ is the radius of the circumscribed circle.

$\alpha_i$  is the half of the angle between the vectors
$\overrightarrow{Op_i}$ and $\overrightarrow{Op}_{i+1}$. The angle
is defined to be positive, orientation is not involved.

$\omega_P$ is the winding number of $P$ with respect to the center
$O$.

$\mu(P)$ is the Morse index of the function $A$ in the point
P. That is, $\mu(P)$ is the number of negative eigenvalues of the
Hessian  $Hess_P(A)$.

A cyclic configuration  is called \textit{central} if one of its
edges contains $O$.

For a non-central configuration, let $\varepsilon_i$ be the
orientation of the edge $p_ip_{i+1}$, that is,

 $\varepsilon_i=\left\{
                       \begin{array}{ll}
                         1, & \hbox{if the center $O$ lies to the left of } p_ip_{i+1};\\
                         -1, & \hbox{if the center $O$ lies to the right of } p_ip_{i+1}.
                       \end{array}
                     \right.$

$E(P)=(\varepsilon_1,\dots,\varepsilon_n)$ is the string of
orientations of all the edges.

$e(P)$ is the number of positive entries in $E(P)$.

\begin{Ex} \label{Expentagon} \cite{panzh} An equilateral pentagonal linkage $L=(1,1,1,1,1)$ has 14 cyclic
configurations  indicated in  Fig.  4.

(1). The convex regular pentagon and its mirror image  are the global maximum and
minimum of the signed area $A$. Their Morse indices are $2$ and $0$
respectively.

(2). The starlike configurations are a local maximum and a local
minimum of $A$.

(3).  There are $10$ more configurations that have three consecutive
edges aligned. Their  Morse indices equal $1$.

\end{Ex}

\begin{figure}[h]\label{pentagon}
\centering
\includegraphics[width=6 cm]{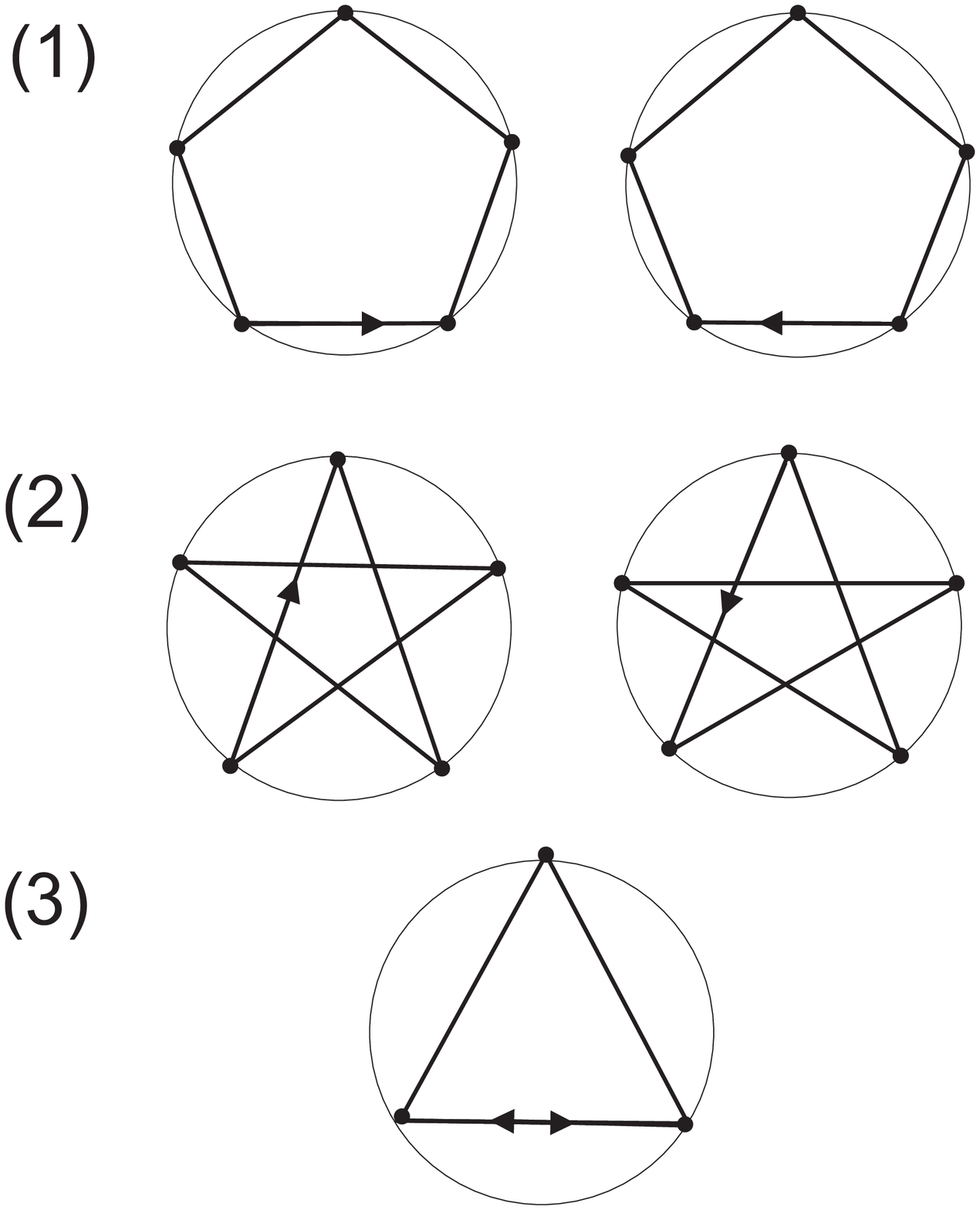}
\caption{Cyclic configurations of the equilateral pentagonal linkage
}
\end{figure}



\section{The decorated moduli space $\widetilde{\mathcal{M}}_3(L)$ and the function $S$}\label{Section_decorated}

We have already defined the moduli space $M_3(L)$. However, it is
convenient to consider the decorated moduli space:

\begin{dfn}
The \textit{decorated moduli space} is defined as the set of pairs
$$\widetilde{M}_3(L)=\{(P,\xi)| P \ \hbox {is a polygon in $\mathbb{R}^3$
 with the sidelengths} \ l_1,...,l_n; \ \xi \in S^2,
\}$$  factorized by the diagonal action of the orientation
preserving isometries of $\mathbb{R}^3$.

Here $S^2 \in \mathbb{R}^3$ is the unit sphere centered at the
origin $O$.
\end{dfn}

\begin{lemma}
The space $\widetilde{M}_3(L)$ is an orientable fibration over $M_3(L)$ whose
fiber is $S^2$.
\end{lemma}
Proof. The set of all polygons with fixed sidelengths (before
factorization by isometries) is known to be orientable. Therefore the set of the
pairs (a polygon, a vector) is also orientable as a trivial
fibration. Since we take a factor by the action
of orientation preserving isometries, the result is also orientable.\qed

\begin{lemma}
The Euler class of the fibration equals zero.
\end{lemma}
Proof. Indeed,  $\xi(P)=\frac{\overrightarrow{p_1p_2}}{|p_1p_2|}$ defines  an everywhere non-zero section. \qed

\begin{cor}\label{Cor_Gisin}(The Gisin sequence for the decorated moduli space)
We have the following short exact sequence:
$$0\rightarrow H^m(M(L)) \rightarrow H^m(\widetilde{M}(L)) \rightarrow H^{m-2}(M(L)) \rightarrow 0.$$
\end{cor}
Proof. This follows directly from Gisin sequence, see \cite{switzer}.
\qed

\begin{dfn} \label{Dfn_areaR3}  Let  $(P,\xi) \in \widetilde{M}_3(L)$, let  $ (x_i,
y_i)$ be the vertices of $P$. The \textit{vector area} of the pair $(P,\xi)$ is defined as
the following scalar product:
$$2S(P,\xi)=( p_1 \times p_2 + p_2 \times p_3+ \dots + p_n\times p_1,\xi).$$

An alternative equivalent definition is:
$$S(P,\xi)=A(pr_{\xi^\perp}(P)),$$
where $pr_{\xi^\perp}$ is the plane orthogonal to $\xi$ and
cooriented by $\xi$.

\end{dfn}

\section{Swap action }

We assume that a polygonal linkage $L$  with all
$l_i$ different is fixed. We make a convention that the numbering is
modulo $n$, that is, for instance, $n+1=1$.

\begin{figure}[h]
\centering
\includegraphics[width=12 cm]{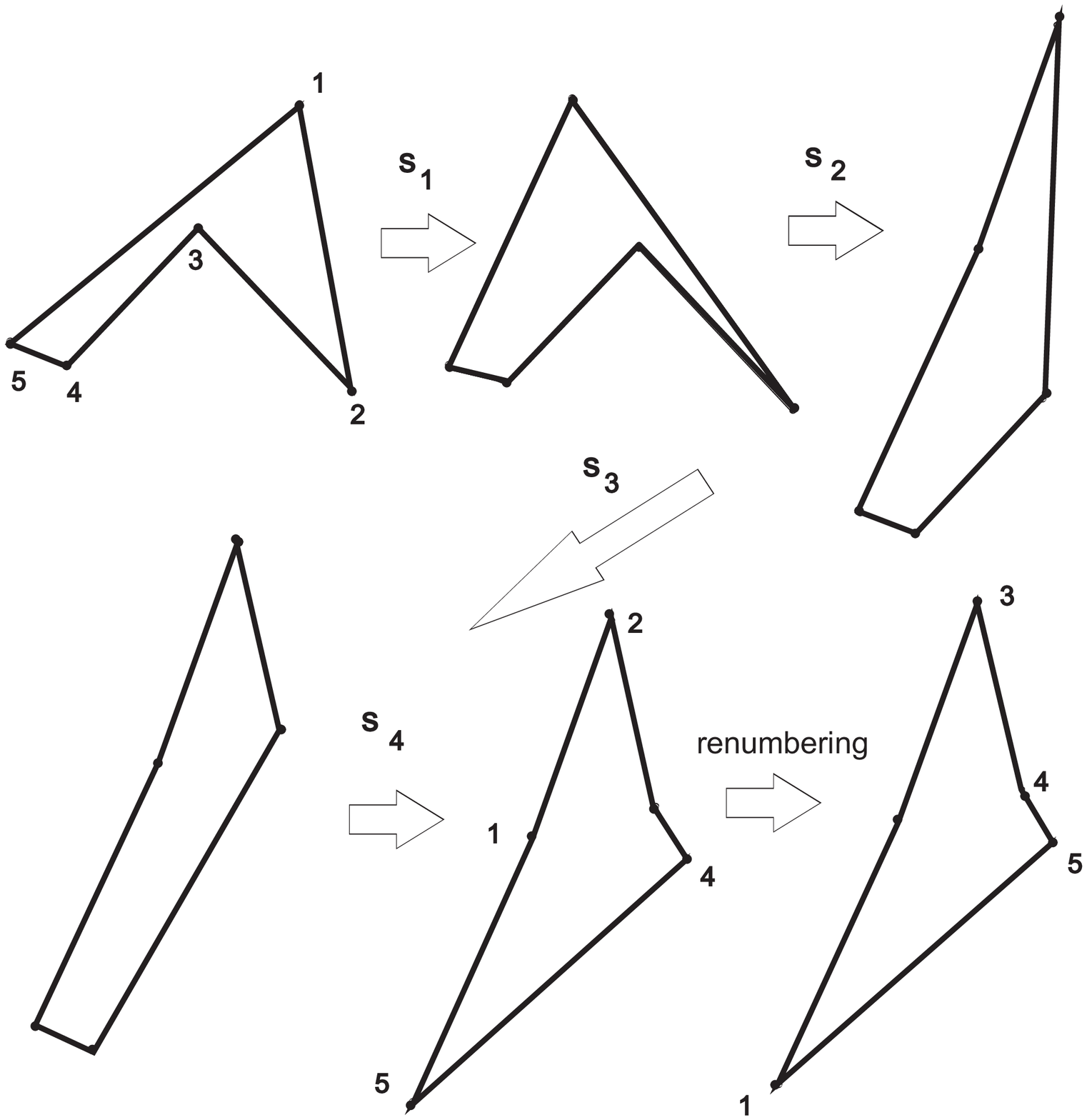}
\caption{The mapping $SW$}\label{Figure_swap}
\end{figure}

\begin{dfn}
\begin{enumerate}
    \item Let $P\in M_{2}(L)$ be a polygon. For $i=1, \dots, n$, denote by
$s_i(P)$ the polygon obtained from $P$ by transposing of the two
edges adjacent to the vertex $p_i$ (see Fig. \ref{Figure_swap}).
    \item For $P\in M_{3}(L)$, the polygon $s_i(P)$ is obtained from
$P$ by the above rules. We assume
that the new pair of edges lies in the plane spanned by the two old
edges.
    \item For $(P, \xi) \in \widetilde{M}_{3}(L)$ we define $s_i(P, \xi)=(s_i(P), \xi).$
\end{enumerate}
\end{dfn}

We get  homeomorphisms $s_i:M_{2,3}(L)\rightarrow M_{2,3}(\sigma_i
L), $ and $s_i:\widetilde{M}_3(L)\rightarrow
\widetilde{M}_3(\sigma_i L), $ where $\sigma_i$ is the element of
the symmetric group $\sigma_i \in S_n$ is a transposition induced by
$s_i$.

Denote also by $SW(P)$  the polygon $s_{n-1}\circ ... \circ s_2
\circ s_1(P)$ whose vertices are renumbered in such a way (that is,
with a shift by one)  that  $SW$ becomes  a (smooth) automorphism of
$M_{2,3}(L)$.

\begin{lemma}
      The actions of $s_i$ and of $SW$ respect the  functions $A$ and $S$.
    \qed
\end{lemma}

\begin{thm}\label{TheoremPlanarSwapInv}
     A polygon $P\in M_{2}(L)$ is $SW$-invariant (that is, $SW(P)$ equals $P$ up
     to an orientation preserving isometry) iff $P$ is cyclic.
    \qed
\end{thm}

\newpage

\section{Critical points  and the Morse index}\label{Section_critical_points}

\begin{figure}[h]
\centering
\includegraphics[width=6 cm]{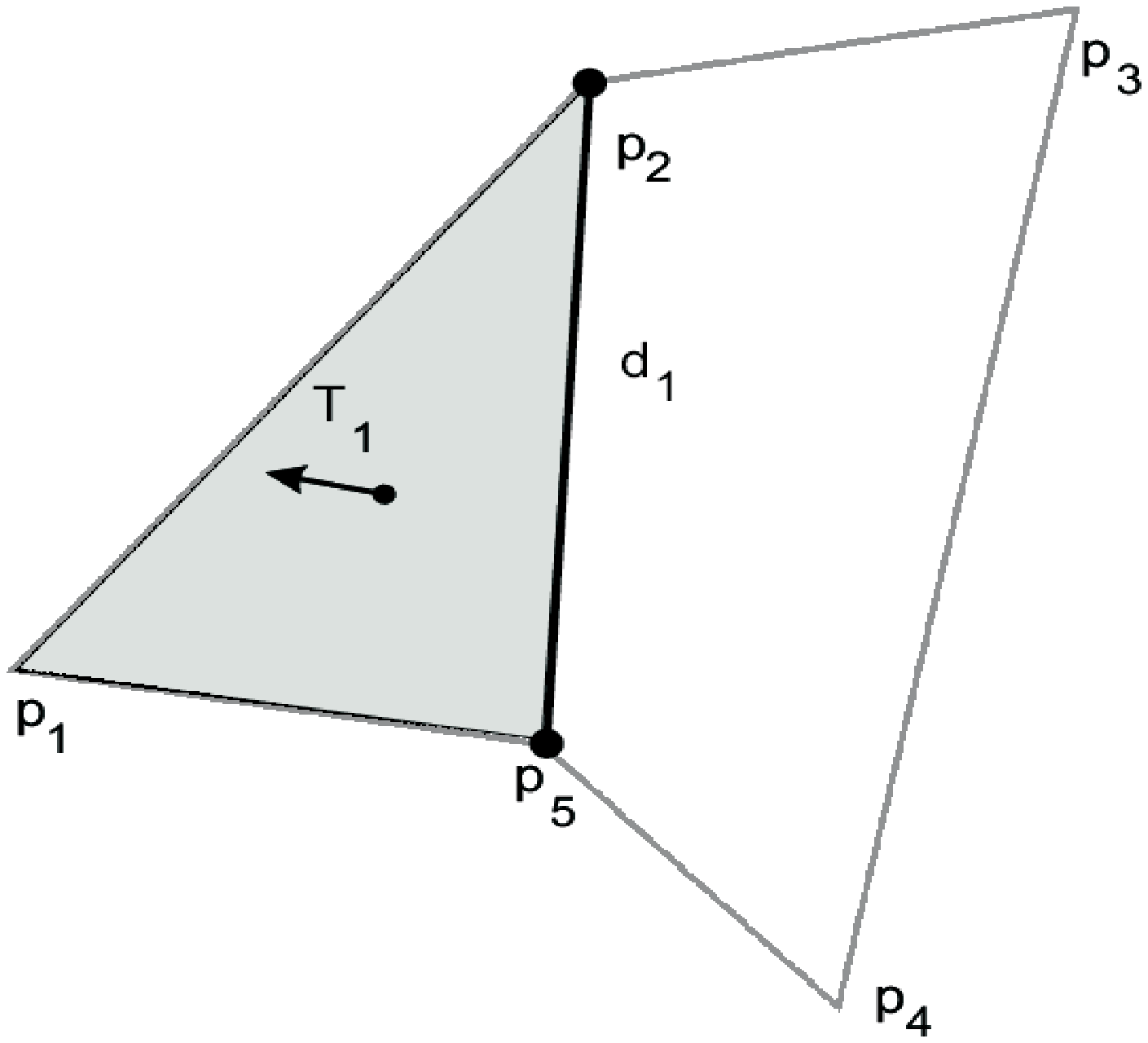}
\caption{Notation for a non-planar critical
polygon}\label{Figure_3Darea}
\end{figure}

\begin{thm}\label{Thm_crirical_3D}
 Generically, critical points $(P, \xi)$ of the function $S$ fall into three classes:
 \begin{itemize}
    \item \textbf{Planar cyclic configurations.}
    These are pairs $(P, \xi)$ such that $P$ is a planar cyclic
    polygon, and $\xi$ is orthogonal to the affine hull of $P$.
    \item \textbf{Non-planar configurations.} They are characterized
    by the three following conditions:
    \begin{enumerate}
\item The vectors $\xi$ and $\overrightarrow{S}=p_1 \times p_2 + p_2 \times p_3+ \dots + p_n\times p_1$
are parallel (but they can have opposite directions).
    \item The orthogonal projection of $P$ onto the plane $\overrightarrow{S(P)}^\perp$ is a
    cyclic polygon.
    \item For every $i$, the vectors $\overrightarrow{T_i}$, $\overrightarrow{S}$, and $\overrightarrow{d_i}$ are
    coplanar.
\end{enumerate}
Here $\overrightarrow{d_i}$ is the $i$-th short diagonal,
$\overrightarrow{T_i}$ is the vector area of the triangle
$p_{i-1}p_ip_{i+1}$, see Fig. \ref{Figure_3Darea}.
    \item \textbf{Zig-zag planar configurations} (existing only for even $n$).
The polygon lies in a plane. There are two parallel lines $l_1$ and
$l_2$ such that all the vertices with even indices lie on the line
$l_1$, whereas all the vertices with odd indices lie on the line
$l_2$. The vector $\xi$ is parallel to $l_{1,2}$.
 \end{itemize}
For all three cases, if $(P, \xi)$ is a critical point, then $(P, -\xi)$ is critical as well.

\end{thm}
Proof.  This follows from \cite{KhristofPanina}, where we proved
nearly the same theorem. \qed

\begin{figure}[h]
\centering
\includegraphics[width=8 cm]{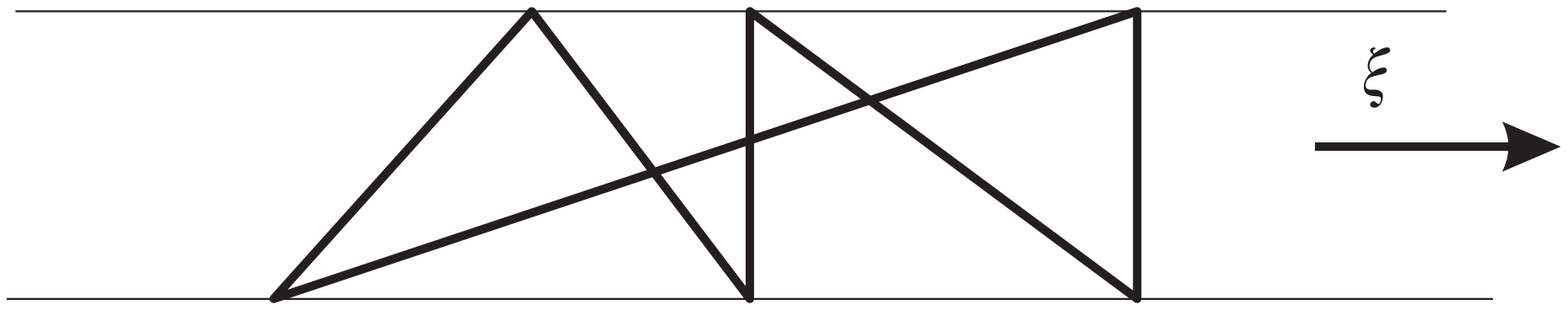}
\caption{A zigzag critical polygon}\label{Figure_zigzag}
\end{figure}

\begin{lemma}\label{LemmaEquilateralPlanar} An equilateral polygon with odd number of edges has no
non-planar critical configurations.
\end{lemma}

A shorter characterization of critical points is given in the
following analogue of Theorem \ref{TheoremPlanarSwapInv}:

\begin{thm}\label{Thm_crirical_3D_SW}
Generically,  $(P, \xi)$ is a critical point of the vector area
function $S$  if and only if  $P$ is SW-invariant.   \qed\end{thm}

An crucial fact about non-planar configurations is the following:
Let $(P, \xi)$ be a non-planar critical configuration. As Theorem
\ref{Thm_crirical_3D} says, the orthogonal projection $prP$ on
$\xi^\perp$ is cyclic.

\begin{lemma}\label{Lemma_deltaZero}
For a non-planar critical configuration,
$$\delta(prP)=0. $$
\end{lemma}
Proof. We use notations $r,pr_i,h_i$ for the radius  of the circle,
lengths of the projections, and  heights differences. We have the
following closing condition $\sum_{i=1}^nh_i=0$.  Note that the
angle $\alpha_i(pr P)$ equals the angle between the chord $pr_i$ and
the circumscribed circle. The conditions from Theorem
\ref{Thm_crirical_3D} imply that the fraction $\frac{h_i}{pr_i}cos \
\alpha _i$ does not depend on $i$. Denote the latter by $h$. The
closing condition implies
$$0=\sum_{i=1}^nh_i=\frac{2h}{r}\sum_{i=1}^n\frac{pr_i/2r}{cos \
\alpha_i}= \frac{2h}{r}\sum_{i=1}^n tan \  \alpha_i.\qed$$

\begin{thm}\label{Thm_Morse_planar}
Let $(P,\xi)$ be a \textbf{planar cyclic} critical point of $S$.

 For the Morse index of the function $S$, we have:
$$\mu(P,\xi)=
  2e(P)-2\omega_P -2  .\qed$$

\end{thm}
The theorem will be proven in the next section.

\bigskip

On the one hand, we
can say nothing about the Morse index of a non-planar critical
polygon. On the other hand, in many cases this result is sufficient
for a construction of a complete Morse theory on the configuration
space. For instance, this is the case for an equilateral polygon
with odd number of edges, see Theorem \ref{LemmaEquilateralTrue}.

\begin{cor}\label{Cor_perfectMorse} If  all SW-invariant
configurations of a polygonal linkage are planar, then
\begin{enumerate}
    \item The function $S$ is a perfect Morse function.
    \item The odd-dimensional homology groups of $M_3(L)$ vanish.
    \item The even-dimensional homology groups are free abelian,
    whose rank can be expressed in terms of the number of cyclic
    configurations of $L$.
\end{enumerate}
\end{cor}
\qed

\section{Proofs for the Section \ref{Section_critical_points}. The equilateral polygon}\label{Section_Proofs}

The Betti numbers and the Euler characteristic of the space $M_3(L)$
are already known due to A. Klyachko. Namely, he proved the
following:
 \begin{thm}\label{ThmKlya1} \cite{klya} The following formulae for the
 Betti numbers are valid:
\begin{eqnarray*} \beta^{2p}( M_3(L))-\beta^{2(p-1)}( M_3(L))&=&\left(\!\!\!\begin{array}{c}n-1\\p\end{array}\!\!\!\right)-\sharp\{I\;
|\; l_I>l/2; |I|=p+1\}=\\&=& \sharp\{I\;|\; l_I<l/2;\;
 |I|=p+1\}-\left(\!\!\!\begin{array}{c}n-1\\p+1\end{array}\!\!\!\right),\end{eqnarray*}

 Here $l=l_1+l_2+\cdots +l_n$, $\ l_I=\sum_{i\in I}l_i$.
 \end{thm}

In the case of equal lengths  the formulae may be simplified.

 \begin{prop} \label{PropKlya2} \cite{klya} For odd $n=2k+1\geq 3$ the Betti numbers of $M_3(1,...,1)$
 are given by the formula

$$\beta^{2p}(M_3(1,...,1)) = \sum _{0\leq i\leq
p}\left(\!\! \begin{array}{c}2k\\ i
\end{array}\!\!\right),\;\;p<k;$$
\end{prop}

\begin{cor}
For Betti numbers of the decorated moduli space
$\widetilde{\beta}^{2p}_n=\beta^{2p}(\widetilde{M}_3(1,...,1))$ we
have
$$\widetilde{\beta}^{2p}_n=\sum _{0\leq i\leq p}\left(\!\!
\begin{array}{c}n\\ i
\end{array}\!\!\right),\;\;p<k;$$
\end{cor}

These expressions can be interpreted
 them  as the
numbers of cyclic equilateral polygons:

\begin{lemma}Let  $n=2k+1$ be an even number.
\begin{enumerate}
    \item  For $p=0,1,..., k$ denote by $N^p_n$ the number of such cyclic
equilateral polygons for which
$$2e-2\omega-2=2p.$$

Then

$$\widetilde{\beta}^{2p}_n=N^p_n.$$
    \item $S$ is a perfect Morse function on the configuration space
$\widetilde{M}_3(1,1,...,1)$.
\end{enumerate}
\end{lemma}

Proof. (1) Indeed, it is easy to find all cyclic equilateral
polygons: first we choose $i$ negatively oriented edges and then
choose the winding number that ranges from $ \pm 1$ to $\pm (k-i)$.

As an illustration, figures \ref{Figure_5gon} and \ref{Figure_7gon}
list all cyclic equilateral pentagons and heptagons.

(2)  By Lemma \ref{Lemma_deltaZero}, equilateral (or, for the sake
of generity, nearly equilateral) polygonal linkages with odd number
of edges have only planar critical configurations. (1) implies
that the number of all critical points equals the sum of all Betti
numbers.\qed

\bigskip

\begin{figure}[h]
\centering
\includegraphics[width=8 cm]{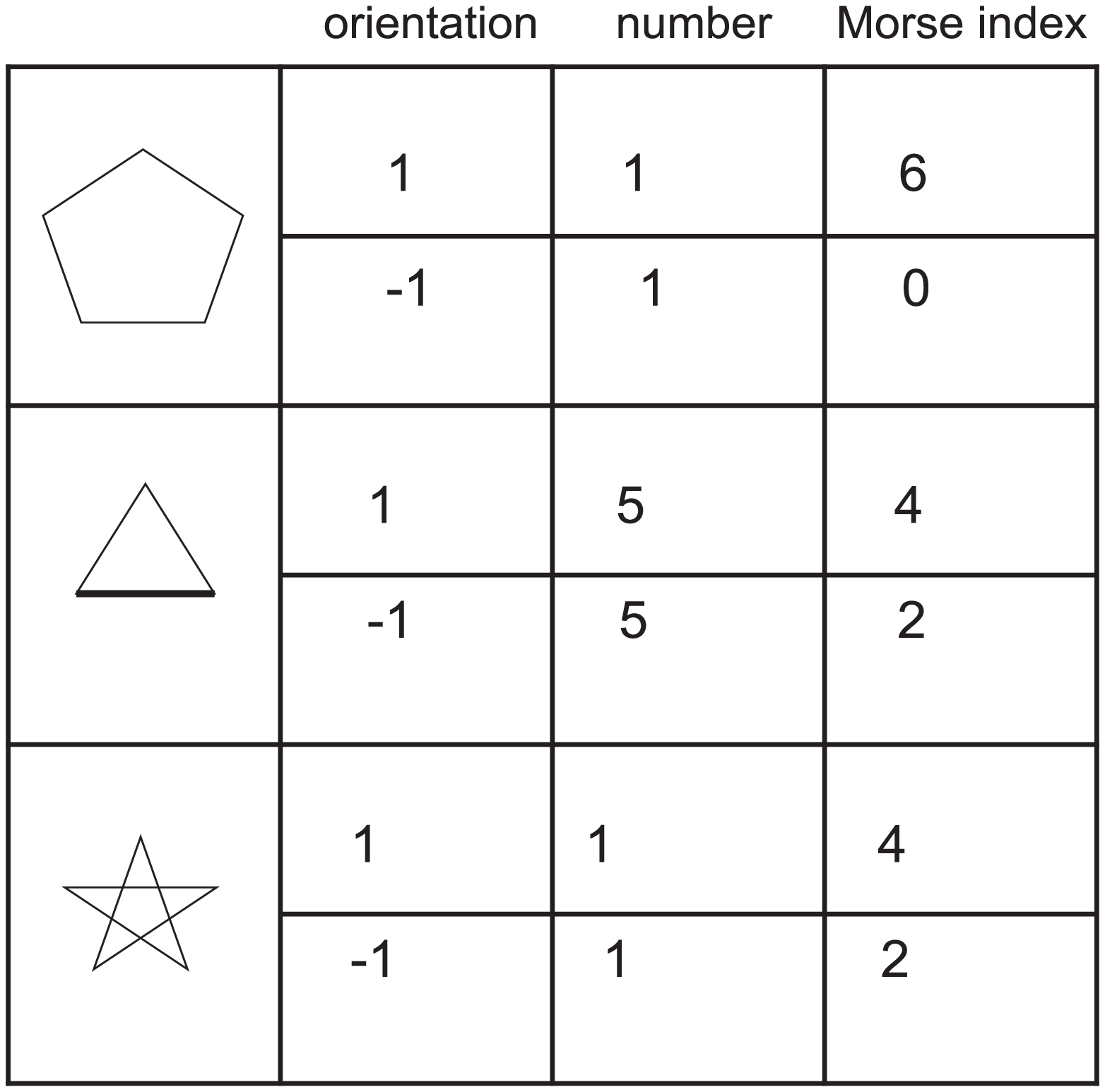}
\caption{Critical equilateral pentagons and their Morse indices}
\end{figure}\label{Figure_5gon}

 \subsection*{Cyclic deformations }
The key idea of how to find the Morse index of a critical point is
to deform the linkage in such a way that the Morse index does not
change.

We start with a planar cyclic configuration $(P,\xi)$. A
\textit{cyclic deformations of a linkage} is a one-parametric
continuous family     arising through the
following construction. Let $(P,\xi)$ be a planar cyclic
configuration of  $L$. We fix the radius $r$ of the circumscribed
circle and the vector $\xi$, and force the vertices $p_i$ to move
along the circle. This yields  a
 continuous family of linkages $L(t)$  together with a
 continuous family of their cyclic configurations $P(t)$.

 We have to understand how the Morse index
$\mu(P,\xi)$ changes during the deformation.

There are only two types of events when $\mu(P,\xi)$
can change:
\begin{enumerate}
    \item If two consecutive vertices $p_i$ and $p_{i+1}$   meet, and the edge
$l_i$ vanishes. This will be called \textit{contraction of the edge
$l_i$}. At such a point the dimension of the configuration space
decreases.
    \item If the point $(P,\xi)$ meets another critical point. If this happens, the value of $\delta(P)$ becomes zero.
\end{enumerate}

\begin{lemma}
If a cyclic deformation $P(t)$ does not pass  through a zero of the
function $\delta$ and has no edge contraction,   the Morse index
$m(P(t))$ remains constant.\qed
\end{lemma}

On the one hand, a detailed analysis of how the Morse index changes
when passing through a zero of $\delta$ provides a proof of Theorem
\ref{Thm_Morse_planar}. This proof is independent on the Klyachko's
result \ref{PropKlya2}.

On the other hand, there exists a shorter proof (the one presented
below), which relies however on the Klyachko's Theorem
\ref{PropKlya2}.

\begin{lemma}\begin{enumerate}
    \item Contraction of  a negatively oriented edge does not change the Morse
index.
    \item Contraction of  a positively oriented edge turns $\mu$ to
$\mu-2$.\qed
\end{enumerate}
\end{lemma}

\begin{lemma}\label{LemmaCyclicDeform}
Let $P=P(0)$ be a planar cyclic polygon. There exists its cyclic
deformation $P(t)$ such that\begin{enumerate}
    \item $\delta (P(t))$ is never zero.
    \item $P(1)$ is an equilateral star with odd number $n=2k+1$ of edges and with $\omega = k$.

For such a deformation, we have
$$\mu(P(1), \xi)=\mu(P, \xi)-2\sharp(\hbox{number of  of positively oriented edges contracted.})$$
\end{enumerate}
Figures 7 and 8  present  examples of such deformations.\qed
\end{lemma}

\begin{figure}[h]
\centering \includegraphics[width=6 cm]{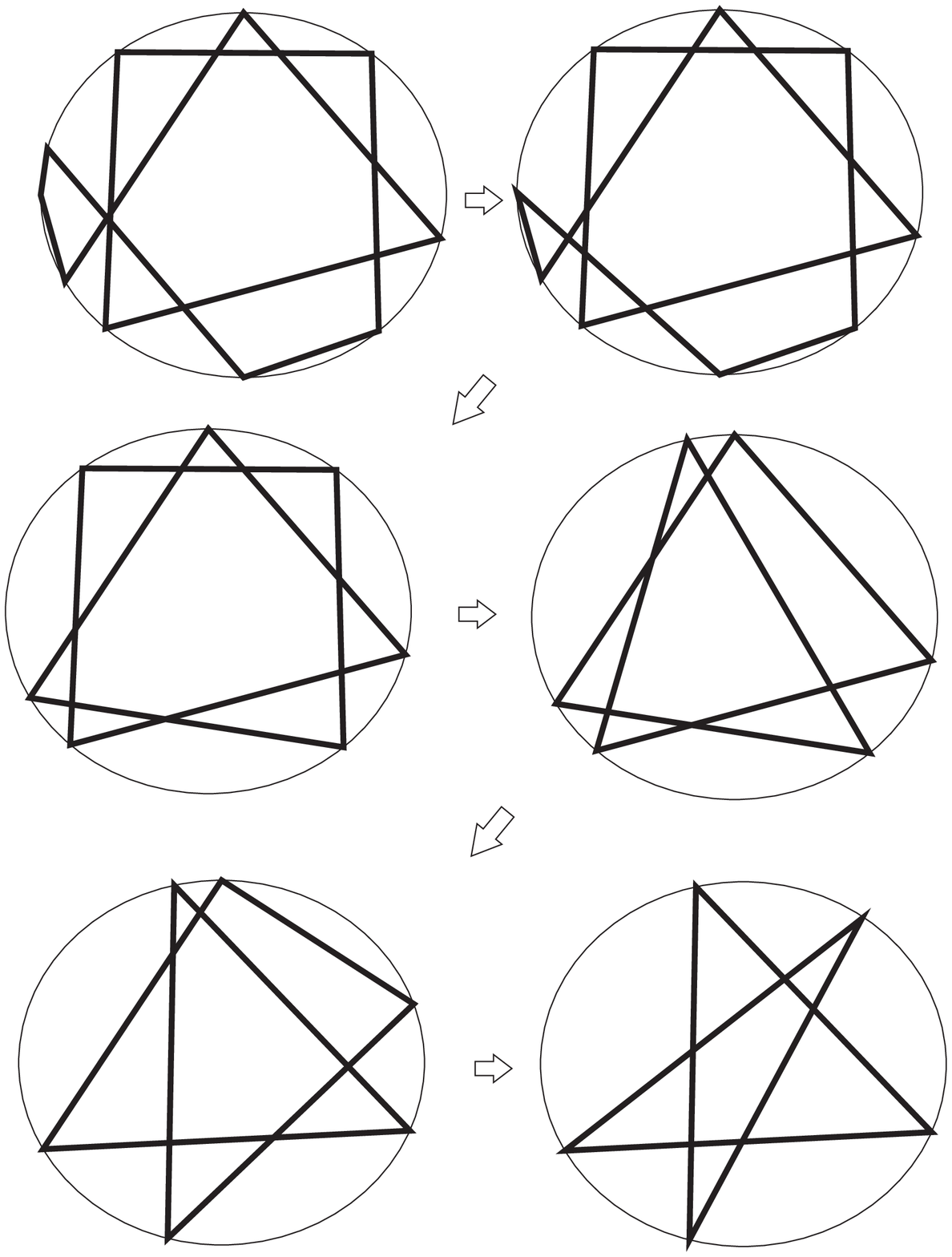}
\caption{A deformation taking a polygon to an equilateral star}\label{FigDeform1}
\end{figure}

\begin{figure}[h]
\centering \includegraphics[width=6 cm]{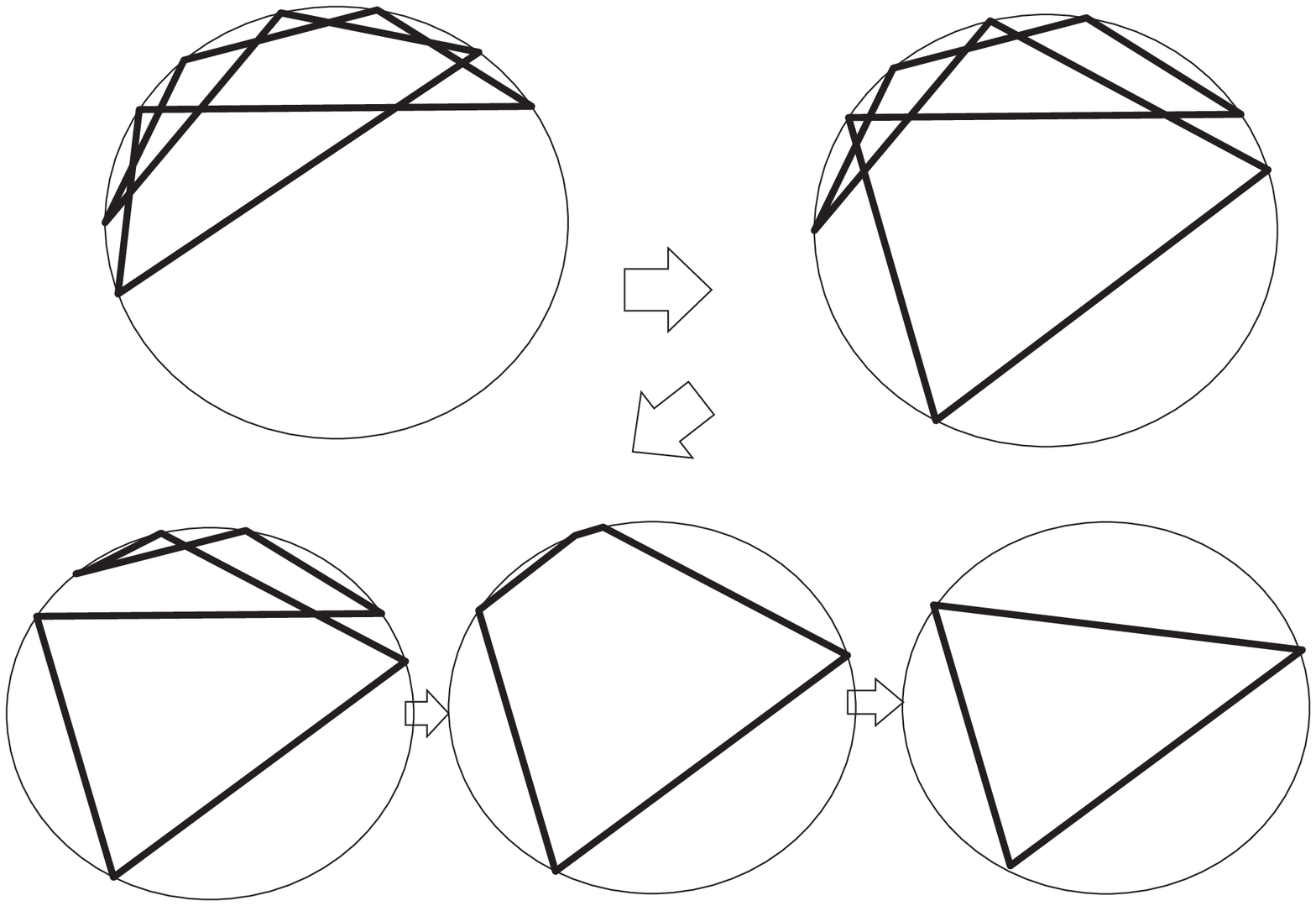}
\caption{One more deformation taking a polygon to a triangle}\label{FigDeform2}
\end{figure}

\begin{thm}\label{LemmaEquilateralTrue}
\begin{enumerate} Let $L=(1,...,1)$ be an equilateral polygons with odd number of
edges. We have the following:
    \item The function $S$ is a perfect Morse function on the decorated moduli space $\widetilde{M}_3(L)$.
    \item The formula for the Morse index from Theorem \ref{Thm_Morse_planar}
is valid for all critical configurations of $L$.
 \item All the Morse indices are even, and the boundary
homomorphisms for the Morse chain complex are zero.
\item The equilateral cyclic  polygons can be interpreted as
independent generators of the homology groups of $\widetilde{M}_3(L)$.
\end{enumerate}
\end{thm}

Proof.  (1) Indeed, by Lemma \ref{LemmaEquilateralPlanar} all
critical configurations are planar. The number of all cyclic
equilateral polygons equals the sum of Betti numbers of the space
$\widetilde{M}_3(L)$.

(2) We prove this inductively by the number of edges. For $n=5$,
this is true by simple reasons.  For induction step, assume that the
statement is proven for $n=2k+1$. Prove it for $n=2k+3$. Lemma
\ref{LemmaCyclicDeform}  determines the Morse indices for the
majority of  the polygons. For instance, the heptagon number  3
(Fig. \ref{Figure_7gon}) after contraction of one negative and one
positive edges gives a pentagonal star whose Morse index is already
known. There are just two polygons that are irreducible in this
sense: the two stars with $\omega =\pm(k+1)$. Since the Morse index
of the positively oriented star is bigger than the Morse index of
the negatively oriented star, the Morse indices are determined
uniquely. The statements (3) and (4) directly follow from (2). \qed


\bigskip

Now we are ready to prove Theorem \ref{Thm_Morse_planar}. Given a
critical point $(P, \xi)$, apply a deformation from Lemma
\ref{LemmaCyclicDeform}. It is easy to check that the difference
$\mu(P(t),\xi)-2e(P(t),\xi)-2\omega (P(t),\xi)-2$ does not change
during the deformation. Besides, by Lemma
\ref{LemmaEquilateralTrue}, the difference is zero for the endpoint
of the deformation. Therefore it is zero at the starting point, that
is, for $(P(t),\xi))$.

\section{More examples}\label{SectionSW_Moreexamples}
\subsection*{An equilateral 7-gon}
\begin{figure}[h]
\centering \includegraphics[width=8.6 cm]{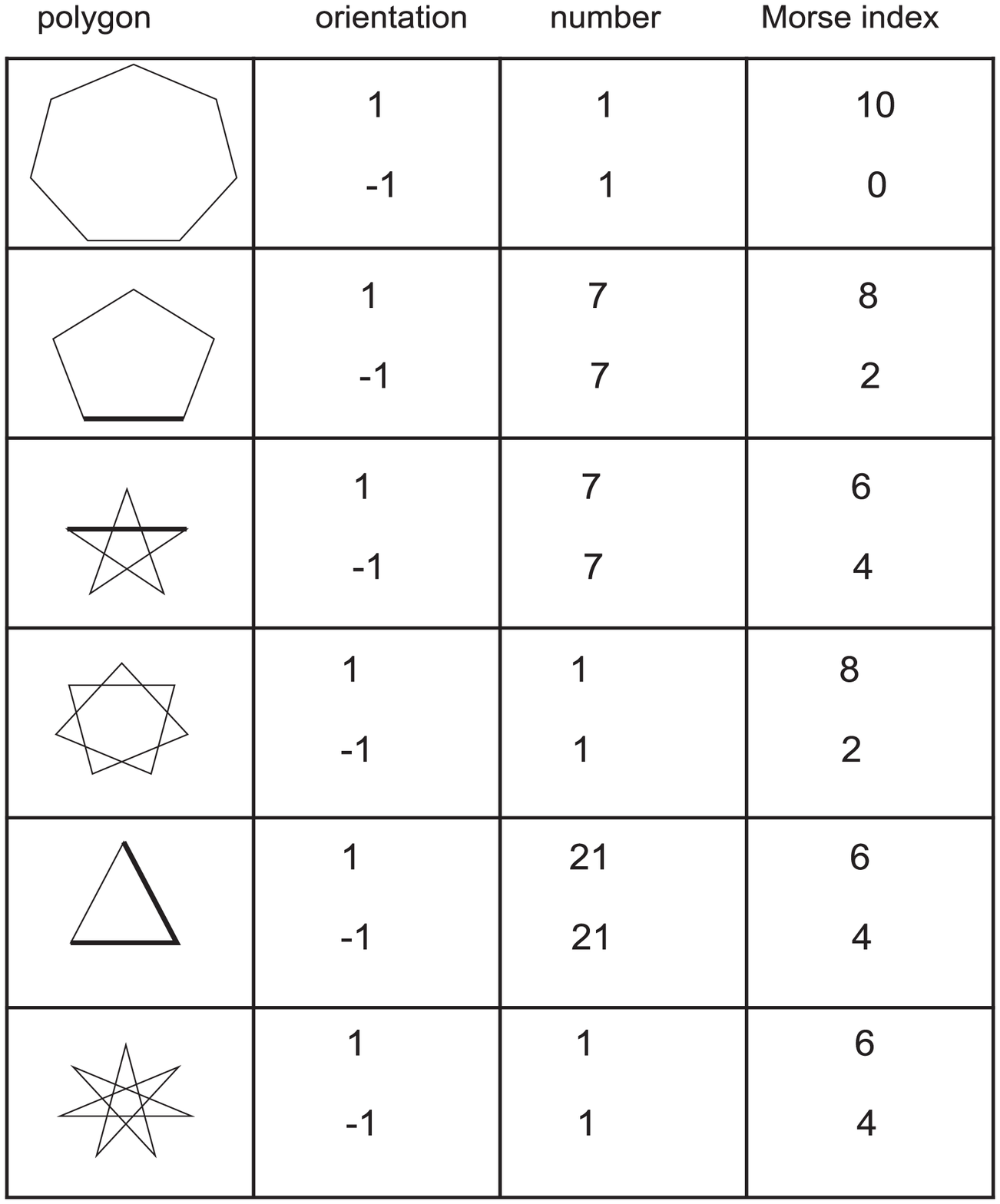}
\caption{Critical equilateral heptagons and their Morse
indices}\label{Figure_7gon}
\end{figure}

Let $L=(1,1,1,1,1,1,1)$ be an equilateral heptagonal linkage. By
Lemma \ref{LemmaEquilateralTrue}, $S$ is a perfect Morse function.
Figure \ref{Figure_7gon}  lists all the types of its critical configurations and their
Morse indices.

\subsection*{A nearly equilateral 6-gon} Let
$L=(1,1,1,1,1,1-\varepsilon)$. Again, $S$ is a perfect Morse
function. Figure \ref{6gon}  lists all the types of its critical configurations
and their Morse indices.

\begin{figure}[h]
\centering
\includegraphics[width=7.5 cm]{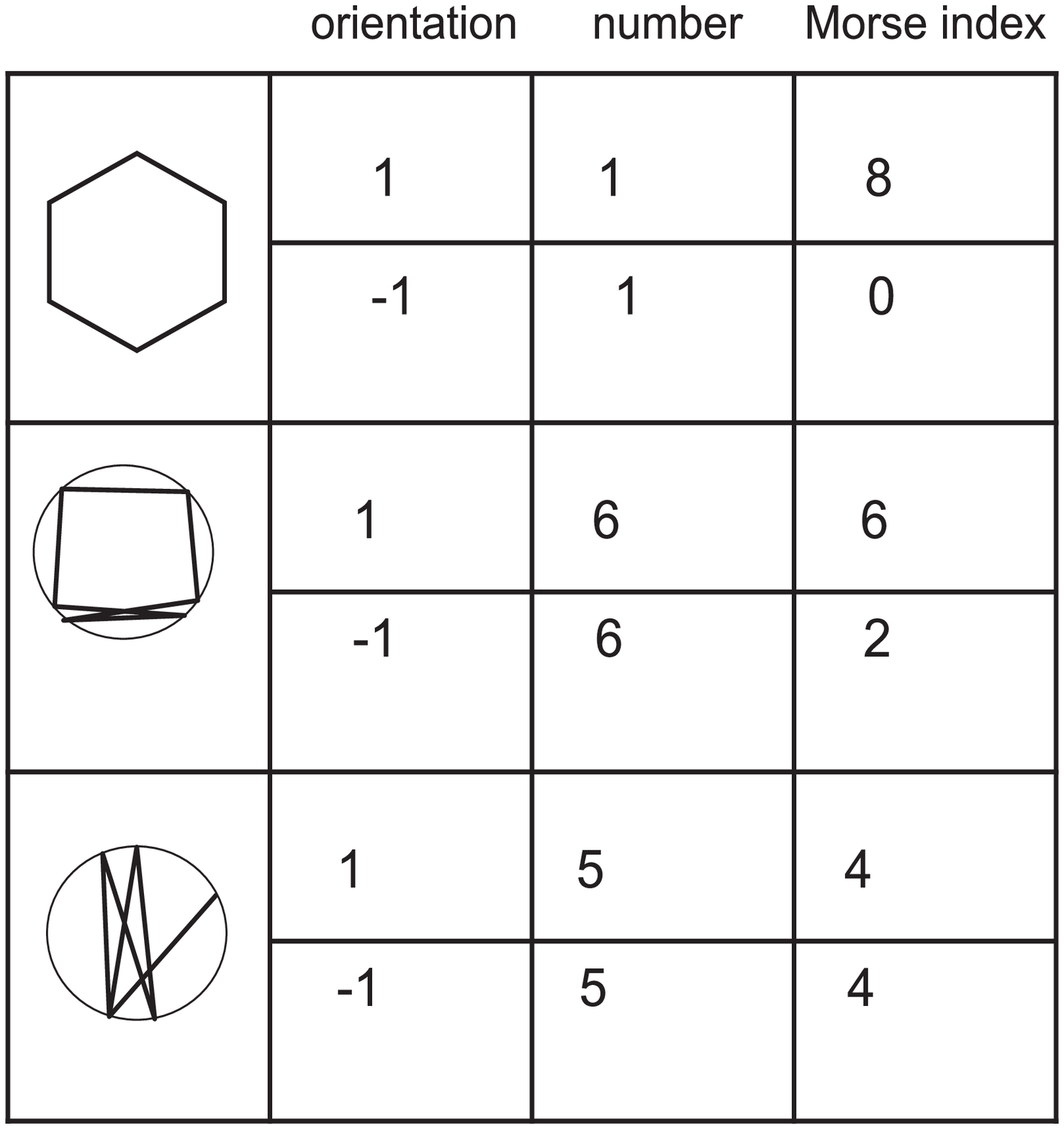}
\caption{Critical nearly equilateral 6-gons and their Morse indices}
\end{figure}\label{6gon}

\newpage

\subsection*{A 4-gonal linkage}
 Let $L=(l_1,l_2,l_3,l_4)$  \ be a generic 4-gonal linkage.
It is known (see \cite{klya}) that $M_3(L)=S^2$. Corollary
\ref{Cor_Gisin}  implies that $H_0(\widetilde{M}_3(L))=
H_4(\widetilde{M}_3(L))=Z, \ H_1(\widetilde{M}_3(L))=
H_3(\widetilde{M}_3(L))=0, \ H_2(\widetilde{M}_3(L))= Z^2$. Let us
establish this result making use of the Morse complex on the space
$\widetilde{M}_3(L)$. There are two possible cases:
\begin{enumerate}
    \item $L$ has only one cyclic configuration which is convex.
    It gives two Morse points with two opposite vectors $\xi$.
    One of them is the global maximum of $S$, and the other one --
    the global minimum. So we have the Morse indices $0$ and $4$.
    Besides, there are exactly two zig-zag critical points with one and the same polygon $P$ and with two opposite vectors
    $\xi$. By symmetry reasons, their Morse indices equal 2.
    \item $L$ has two cyclic configuration, one is convex and the other one is self-intersecting.
    Each of them gives two Morse points with two opposite vectors $\xi$.
    The convex polygon yields the global maximum and the
    the global minimum. As in the previous case, we have the Morse indices $0$ and $4$.
    The self-intersecting configuration yields two critical points with Morse indices equal 2.
\end{enumerate}

In both cases the Morse chain complex has zero chain groups with odd
indices. Therefore, The odd homology groups are zero, and the even
homology groups are free abelian whose rank equals the number of
critical points.

\subsection*{A 5-gonal linkage with just two planar configurations}
Consider a polygonal linkage $$L=(1,1,1,1,4-\varepsilon)$$ where
$\varepsilon$ is small. There are  four critical points: two planar
ones (the global maximum and the global minimum of $S$) and two
non-planar ones that differ on a mirror symmetry with respect to a
plane. Again, $S$ is a perfect Morse function.

\end{document}